\documentclass[10pt,reqno]{amsart}

\usepackage{amsmath,amsthm}
\usepackage{amssymb,amsfonts}
\usepackage{euscript}

\numberwithin{equation}{subsection}

\newcommand{\sqsp}{\renewcommand{\baselinestretch}{1.25}\tiny\normalsize}

\newcommand{\relphantom[1]}{\mathrel{\phantom{#1}}}

\newtheorem{thm}[subsection]{Theorem}
\newtheorem{lemma}[subsection]{Lemma}
\newtheorem{prop}[subsection]{Proposition}
\newtheorem{cor}[subsection]{Corollary}

\theoremstyle{definition}

\newcommand{\cat}[1]{{\EuScript #1}}
\newcommand{\cA}{\cat{A}}

\newcommand{\cF}{\cat{F}}

\newcommand{\bF}{\mathbf{F}}

\newcommand{\Lambdabar}{\overline{\Lambda}}
\newcommand{\Lambdatilde}{\widetilde{\Lambda}}
\newcommand{\lambdabar}{\overline{\lambda}}

\newcommand{\sumprime}{\sideset{}{'}\sum}  
\newcommand{\Rbar}{\overline{R}}

\DeclareMathOperator{\Id}{Id}
\DeclareMathOperator{\Hom}{Hom}
\DeclareMathOperator{\End}{End}
\DeclareMathOperator{\Ob}{Ob}
\DeclareMathOperator{\Der}{Der}
\DeclareMathOperator{\Coder}{Coder}

\begin{document}
\title{Cohomology and deformation of module-algebras}
\author{Donald Yau}

\begin{abstract}
An algebraic deformation theory of module-algebras over a bialgebra is constructed.  The cases of module-coalgebras, comodule-algebras, and comodule-coalgebras are also considered.
\end{abstract}

\address{Department of Mathematics, The Ohio State University Newark, 1179 University Drive, Newark, OH 43055, USA}
\email{dyau@math.ohio-state.edu}
\maketitle
\sqsp


\section{Introduction}
\label{sec:intro}


Let $H$ be a bialgebra and $A$ be an associative algebra.  The algebra $A$ is said to be an $H$-module-algebra if there is an $H$-module structure on $A$ such that the multiplication on $A$ becomes an $H$-module morphism.  There are many important examples of this structure.  For example, an algebra over the Landweber-Novikov algebra $S$ \cite{landweber,novikov} is an $S$-module-algebra.  In particular, the complex cobordism $MU^*(X)$ of a topological space $X$, equipped with its commutative ring structure and the stable cobordism operations, is an $S$-module-algebra.  Likewise, if $p$ is a prime, then algebras over the Steenrod algebra $\cA_p$ \cite{es,milnor} are $\cA_p$-module-algebras.  In particular, the singular mod $p$ cohomology $H^*(X; \bF_p)$ of a topological space $X$, equipped with its commutative $\bF_p$-algebra structure and the Steenrod operations, is an $\cA_p$-module-algebra.  Other examples from algebraic topology can be found in \cite{boardman}.  There are also important examples from Lie and Hopf algebras theory.  For instance, the affine plane admits a module-algebra structure over the enveloping bialgebra of the Lie algebra $sl(2)$ \cite[V.6]{kassel}.


The main purposes of this paper are:
\begin{enumerate}
\item Construct the deformation cohomology for any module-algebra $A$ over a bialgebra $H$, where the deformation is taken with respect to the $H$-module structure, keeping the algebra structure on $A$ unaltered.
\item Use the deformation cohomology to describe infinitesimal, rigidity, and extension results concerning algebraic deformations of module-algebras.
\end{enumerate}
We will also consider the related cases of module-coalgebras, comodule-algebras, and comodule-coalgebras.  These three other algebraic structures are important in the studies of Hopf algebras and quantum groups (see, e.g., \cite{d,kassel,mont,sweedler}).


Some remarks are in order.  In \cite{yau} the author considered deformations of algebras over the Landweber-Novikov algebra, which, as mentioned above, are examples of module-algebras.  The current paper generalizes the constructions and results there to any module-algebra.  Moreover, in the study of the obstructions to extending cocycles to deformations, we actually obtain a simpler and more conceptual argument.  More precisely, in \cite[Lemma 5.2]{yau}, a certain obstruction class was shown to be a $1$-cocycle in the deformation complex using a rather computational argument. It makes heavy use of the composition law and the Cartan formula for the Landweber-Novikov operations.  In the current paper, the corresponding fact is proved by a much simpler argument, using a certain cup-product in the deformation complex.  This brings us to the next remark.

In most known cases of algebraic deformations, the cochain complex that controls the deformations is a dg Lie algebra (see, e.g., \cite{bal,ger1,hinich}).  The Lie bracket is usually used when one tries to show that the obstructions to extending cocycles to deformations are themselves cocycles.  On the other hand, the cochain complex that controls the deformations of a module-algebra is a differential graded algebra (DGA), whose product we denote by a cup-product.  Moreover, part of the building blocks for this DGA is a certain Hochschild cochain complex.  The cup-product in the deformation DGA extends the usual cup-product in this Hochschild cochain complex.  It remains an open question as to whether the cohomology of the deformation DGA admits a Gerstenhaber-algebra structure, like the one on Hochschild cohomology \cite{ger0}.

There may be a more general deformation theory of the pair $(H,A)$, in which, say, both the $H$-module structure and the algebra structure on $A$ are deformed.  The resulting deformation complex $\cF^\bullet(H,A)$ should contain the deformation DGA $\cF^\bullet(A)$ in this paper.  The extended deformation complex $\cF^\bullet(H,A)$ may then be a dg Lie algebra.  This idea was pointed out by the referee.

\subsection{Organization}
\label{subsec:org}

The deformation DGA $\cF^\bullet(A)$ for a module-algebra $A$ is constructed in the next section.  In section \ref{sec:deformation}, we use the deformation DGA to study deformations, infinitesimals, and rigidity of a module-algebra.  Infinitesimals are identified with certain cohomology classes in $H^1(\cF^\bullet(A))$ (Theorem \ref{thm:inf}).  The rigidity result (Corollary \ref{cor:rigidity}) states that every deformation is equivalent to the trivial one, provided that both $H^1(\cF^\bullet(A))$ and $H^2_h(A,A)$ are trivial.  Here $H^\bullet_h(A,A)$ is the Hochschild cohomology of $A$ with coefficients in itself.  In section \ref{sec:ext}, the obstructions to extending a $1$-cocycle in the deformation DGA to a full-blown deformation are identified.  They are shown to be $2$-cocycles (Lemma \ref{lem:ob}).  In particular, if $H^2$ of the deformation DGA is trivial, then such a deformation always exists (Corollary \ref{cor:ext}).  In section \ref{sec:def'}, the deformation DGAs for a module-coalgebra, a comodule-algebra, and a comodule-coalgebra are constructed and the corresponding results about deformations are listed.

\subsection{Acknowledgment}
\label{subsec:ack}

The author thanks the referee for reading an earlier version of this paper and for his/her constructive comments and suggestions.


\section{Deformation cohomology for module-algebras}
\label{sec:def coh}

The purpose of this section is to construct the DGA that controls the deformations of a module-algebra over a bialgebra.

\subsection{Notations}
\label{subsec:notations}

Fix a ground field $K$.  Tensor products and $\Hom$ will be taken over $K$. Also fix a $K$-bialgebra $(H, \mu_H, \Delta_H)$.

Denote by $\End(X)$ the algebra, under composition, of $K$-linear endomorphisms of a vector space $X$.  For an algebra $(A, \mu_A)$, a \emph{derivation on $A$} is a linear map $\varphi \in \End(A)$ such that $\varphi \circ \mu_A = \mu_A \circ (\Id_A \otimes \varphi + \varphi \otimes \Id_A)$.  The set of derivations on $A$ is denoted by $\Der(A)$.

In a coalgebra $(C, \Delta)$, we use Sweedler's notation \cite{sweedler} for comultiplication: $\Delta(x) = \sum_{(x)} x_{(1)} \otimes x_{(2)}$, $\Delta^2(x) = \sum_{(x)} x_{(1)} \otimes x_{(2)} \otimes x_{(3)}$, etc.  A \emph{coderivation on $C$} is a linear map $\varphi \in \End(C)$ such that $\Delta \circ \varphi = (\Id_C \otimes \varphi + \varphi \otimes \Id_C) \circ \Delta$.  The set of coderivations on $C$ is denoted by $\Coder(C)$.

\subsection{Module-algebras}
\label{subsec:module-alg}

Basic information about module-algebras can be found in many books on Hopf algebras, e.g., \cite{d,kassel,mont,sweedler}.  Let $(A, \mu_A)$ be an associative $K$-algebra.  Say that $A$ is an $H$-\emph{module-algebra} if and only if there exists an $H$-module structure $\lambda \in \Hom(H, \End(A))$ on $A$ such that $\mu_A \colon A \otimes A \to A$ is an $H$-module morphism.  In other words, $\lambda$ is required to satisfy:
   \[
   \begin{split}
   \lambda(xy) & \,=\, \lambda(x) \circ \lambda(y), \\
   \lambda(x)(ab) & \,=\, \sum_{(x)} \lambda(x_{(1)})(a) \cdot \lambda(x_{(2)})(b),
   \end{split}
   \]
for $x, y \in H$ and $a, b \in A$.

For this and the next two sections, $A$ will be an $H$-module-algebra with $H$-action map $\lambda$.

\subsection{Hochschild cohomology}
\label{subsec:Hoch}

The deformation complex of a module-algebra uses the Hochschild cochain complex \cite{hochschild}, which we now recall.  Let $(R, \mu)$ be an algebra and let $M$ be an $R$-bimodule with left $R$-action $\alpha_L$ and right $R$-action $\alpha_R$.  For integers $n \geq 0$, the \emph{module of Hochschild $n$-cochains} of $R$ with coefficients in $M$ is $C^n_h(R,M) = \Hom(R^{\otimes n}, M)$.  The differential $b \colon C^n_h(R,M) \to C^{n+1}_h(R,M)$ is given by
   \begin{multline*}
   b\varphi \,=\, \alpha_L \circ (\Id_R \otimes \varphi) \,+\, \sum_{i=1}^n (-1)^i \varphi \circ \left( \Id_{R^{\otimes(i-1)}} \otimes \mu \otimes \Id_{R^{\otimes(n-i)}} \right) \\
   +\, (-1)^{n+1} \alpha_R \circ (\varphi \otimes \Id_R)
   \end{multline*}
for $\varphi \in C^n_h(R,M)$.  The \emph{Hochschild cohomology of $R$ with coefficients in $M$} is the cohomology of the cochain complex $C^\bullet_h(R,M)$ and is denoted by $H^\bullet_h(R,M)$.

\subsection{Deformation complex}
\label{subsec:def complex}

Now we define the deformation complex $(\cF^\bullet(A), d)$ of $A$ as an $H$-module-algebra.  Set:
   \[
   \begin{split}
   \cF^0(A) & \,=\, \Der(A), \\
   \cF^1(A) & \,=\, \Hom(H, \End(A)).
   \end{split}
   \]
For integers $n \geq 2$, define
   \[
   \cF^n(A) \,=\, \cF^n_0(A) \oplus \cF^n_1(A),
   \]
where
   \[
   \begin{split}
   \cF^n_0(A) & \,=\, C^n_h(H, \End(A)), \\
   \cF^n_1(A) & \,=\, \Hom(H, \Hom(A^{\otimes n}, A)).
   \end{split}
   \]
In $C^n_h(H, \End(A))$, we consider $\End(A)$ as an $(H, \mu_H)$-bimodule via the structure map
   \[
   \begin{split}
   H \otimes \End(A) \otimes H & \to \End(A) \\
   x \otimes f \otimes y & \mapsto \lambda(x) \circ f \circ \lambda(y).
   \end{split}
   \]

Now we define the differentials $d \colon \cF^n(A) \to \cF^{n+1}(A)$.  For $\varphi \in \Der(A)$, set
   \[
   d^0\varphi \,=\, \lambda \circ \varphi - \varphi \circ \lambda,
   \]
where $(\lambda \circ \varphi)(x) = \lambda(x) \circ \varphi$ and $(\varphi \circ \lambda)(x) = \varphi \circ \lambda(x)$ for $x \in H$.  For integers $n \geq 1$, set
   \[
   d^n \,=\, (d^n_0; d^n_1),
   \]
where $d^n_i \colon \cF^n_i(A) \to \cF^{n+1}_i(A)$ $(i = 0, 1)$ is defined as follows:
   \begin{itemize}
   \item $d^n_0 = b \colon C^n_h(H, \End(A)) \to C^{n+1}_h(H, \End(A))$, the Hochschild coboundary.
   \item $d^n_1 = \sum_{i=0}^{n+1} (-1)^i d^n_1 \lbrack i \rbrack$, where
   \[
   (d^n_1 \lbrack i \rbrack)(\varphi_1)(x)(\mathbf{a}) \,=\,
   \begin{cases}
   \sum_{(x)} \lambda(x_{(1)})(a_1) \cdot \varphi_1(x_{(2)})(a_2 \otimes \cdots \otimes a_{n+1}) & \text{if } i = 0 \\
   \varphi_1(x)(a_1 \otimes \cdots \otimes (a_i a_{i+1}) \otimes \cdots a_{n+1}) & \text{if } 1 \leq i \leq n \\
   \sum_{(x)} \varphi_1(x_{(1)})(a_1 \otimes \cdots \otimes a_n) \cdot \lambda(x_{(2)})(a_{n+1}) & \text{if } i = n + 1. \end{cases}
   \]
   \end{itemize}
Here $\varphi = (\varphi_0; \varphi_1) \in \cF^n(A)$, $x \in H$, and $\mathbf{a} = a_1 \otimes \cdots \otimes a_{n+1} \in A^{\otimes n+1}$.  In these definitions when $n = 1$, we think of $\cF^1_0 = \cF^1_1 = \cF^1$ and $\varphi_0 = \varphi_1 = \varphi$.

\begin{prop}
\label{prop:def complex}
$(\cF^\bullet(A), d)$ is a cochain complex
\end{prop}

\begin{proof}
One can check directly that $d^1 \circ d^0 = 0$.  It is clear that $d^{i+1}_0 \circ d^i_0 = 0$ for $i \geq 1$, since $d^i_0 = b$ is the Hochschild coboundary.  For $0 \leq k < l \leq i + 2$, it is straightforward to check the cosimplicial identities
   \[
   d^{i+1}_1 \lbrack l \rbrack \circ d^i_1 \lbrack k \rbrack \,=\,
   d^{i+1}_1 \lbrack k \rbrack \circ d^i_1 \lbrack l - 1 \rbrack.
   \]
As usual, this implies that $d^{i+1}_1 \circ d^i_1 = 0$.
\end{proof}

\subsection{Cup-product}
\label{subsec:cup}

The usual cup-product on the Hochschild cochain complex $C^\bullet_h(H, \End(A))$ is defined as
   \[
   (f \cup g)(x_1 \otimes \cdots \otimes x_{m+n}) \,=\,
   f(x_1 \otimes \cdots \otimes x_m) \circ g(x_{m+1} \otimes \cdots \otimes x_{m+n})
   \]
for $f \in C^m_h(H, \End(A))$, $g \in C^n_h(H, \End(A))$, and $x_1, \ldots, x_{m+n} \in H$.

Using this, we define a cup-pairing
   \begin{equation}
   \label{eq:cup}
   - \cup - \colon \cF^m(A) \otimes \cF^n(A) \to \cF^{m+n}(A)
   \end{equation}
for integers $m, n > 0$ as follows.  (Note that we do not consider the cases where $m = 0$ or $n = 0$.)  Suppose that $f = (f_0; f_1) \in \cF^m(A)$ and $g =(g_0; g_1) \in \cF^n(A)$ for $m, n > 0$.  Define:
   \begin{itemize}
   \item $(f \cup g)_0 = f_0 \cup g_0$, where the $\cup$-product on the right-hand side is the one on $C^\bullet_h(H, \End(A))$.
   \item For $x \in H$ and $\mathbf{a} = a_1 \otimes \cdots \otimes a_{m+n} \in A^{\otimes m + n}$,
   \[
   (f \cup g)_1(x)(\mathbf{a}) \,=\, \sum_{(x)} f_1(x_{(1)})(a_1 \otimes \cdots \otimes a_m) \cdot g_1(x_{(2)})(a_{m+1} \otimes \cdots \otimes a_{m+n}).
   \]
   \end{itemize}
In these definitions, if $m = 1$, then we think of $f_0 = f_1 = f$, and similarly when $n = 1$.

\begin{prop}
\label{prop:cup}
The $\cup$-pairing in \eqref{eq:cup} is associative and satisfies the Leibniz identity,
   \[
   d(f \cup g) \,=\, (df) \cup g + (-1)^{\vert f \vert} f \cup (dg),
   \]
where $\vert f \vert = m$ for $f \in \cF^m(A)$.  In particular, it follows that $(\cF^{> 0}(A), d, \cup)$ is a DGA.
\end{prop}

\begin{proof}
The associativity of the $\cup$-product on $C^\bullet_h(H, \End(A))$ is obvious from the definition.  The associativity of $\cup$ on the component $\cF^\bullet_1(A)$ is an immediate consequence of the coassociativity of $\Delta_H$, namely,
   \[
   \sum_{(x)(x_{(2)})} x_{(1)} \otimes (x_{(2)})_{(1)} \otimes (x_{(2)})_{(2)}
   = \Delta_H^2(x)
   = \sum_{(x)(x_{(1)})} (x_{(1)})_{(1)} \otimes (x_{(1)})_{(2)} \otimes x_{(2)}.
   \]
The Leibniz identity can be seen by direct inspection.
\end{proof}

This Proposition implies, as usual, that the $\cup$-pairing descends to cohomology.

\begin{cor}
\label{cor:cup}
The $\cup$-pairing \eqref{eq:cup} induces a well-defined product
   \[
   - \cup - \colon H^m(\cF^\bullet(A)) \otimes H^n(\cF^\bullet(A)) \to H^{m+n}(\cF^{\bullet}(A))
   \]
for $m, n > 0$, making $H^{\geq 1}(\cF^\bullet(A))$ into a graded algebra.
\end{cor}



\section{Formal deformation and rigidity}
\label{sec:deformation}

The purposes of this section are to (i) define deformations, (ii) identify infinitesimals with suitable cohomology classes, and (iii) obtain a cohomological criterion for rigidity.

\subsection{Deformation}
\label{subsec:deformation}

By a \emph{deformation of $A$} (as an $H$-module-algebra), we mean a power series $\Lambda_t = \sum_{i=0}^\infty \lambda_i t^i$ with $\lambda_0 = \lambda$ and each $\lambda_i \in \cF^1(A)$, satisfying
   \begin{subequations}
   \label{eq:deformation}
   \begin{align}
   \Lambda_t(xy) & \,=\, \Lambda_t(x) \circ \Lambda_t(y), \label{eq1:def} \\
   \Lambda_t(x)(ab) & \,=\, \sum_{(x)} \Lambda_t(x_{(1)})(a) \cdot \Lambda_t(x_{(2)})(b) \label{eq2:def}
   \end{align}
   \end{subequations}
for $x, y \in H$ and $a, b \in A$.  In particular, by linearity, such a $\Lambda_t$ gives $(A \lbrack \lbrack t \rbrack \rbrack, \mu_A)$ an $H$-module-algebra structure, which reduces to the original one when setting $t = 0$.  The linear term $\lambda_1$ is called the \emph{infinitesimal} of $\Lambda_t$.

In order to identify the infinitesimal with a suitable cohomology class, we need an appropriate notion of equivalence.

\subsection{Equivalence}
\label{subsec:equivalence}

A \emph{formal automorphism of $A$} is a power series $\Phi_t = \sum_{i=0}^\infty \phi_i t^i$ with $\phi_0 = \Id_A$ and each $\phi_i \in \End(A)$ such that $\Phi_t$ is multiplicative, i.e.
   \begin{equation}
   \label{eq:mult}
   \Phi_t(ab) \,=\, \Phi_t(a) \Phi_t(b)
   \end{equation}
for all $a, b \in A$.

Note that this is exactly the same definition as in the special case of algebras over the Landweber-Novikov algebra \cite[3.2]{yau}.

Suppose that $\Lambda_t = \sum_{i=0}^\infty \lambda_i t^i$ and $\Lambdabar_t = \sum_{i=0}^\infty \lambdabar_i t^i$ are deformations of $A$.  We say that $\Lambda_t$ and $\Lambdabar_t$ are \emph{equivalent} if and only if there exists a formal automorphism $\Phi_t$ of $A$ such that
   \begin{equation}
   \label{eq:equiv}
   \Lambdabar_t \,=\, \Phi^{-1}_t \Lambda_t \Phi_t.
   \end{equation}
On the right-hand side, one considers $\phi_i \lambda_j \phi_k$ as an element of $\cF^1(A)$ via the formula,
   \[
   (\phi_i \lambda_j \phi_k)(x) = \phi_i \circ \lambda_j(x) \circ \phi_k
   \]
for $x \in H$.  It is clear that this is an equivalence relation.  Moreover, given a deformation $\Lambda_t$ and a formal automorphism $\Phi_t$, one can define another deformation $\Lambdabar_t$ using \eqref{eq:equiv}, and it is automatically equivalent to $\Lambda_t$.

The following result properly identifies the infinitesimal of a deformation with a cohomology class.

\begin{thm}
\label{thm:inf}
Let $\Lambda_t = \sum_{i=0}^\infty \lambda_i t^i$ be a deformation of $A$.  Then $\lambda_1 \in \cF^1(A)$ is a $1$-cocycle whose cohomology class is determined by the equivalence class of $\Lambda_t$.  Moreover, if $\lambda_1 = \cdots = \lambda_k = 0$, then $\lambda_{k+1}$ is a $1$-cocycle in $\cF^1(A)$.
\end{thm}

\begin{proof}
The condition \eqref{eq1:def} is equivalent to the equality
   \begin{equation}
   \label{eq:def1}
   \lambda_n(xy) \,=\, \sum_{i+j\,=\,n} \lambda_i(x) \circ \lambda_j(y)
   \end{equation}
for all $n \geq 0$ and $x, y \in H$.  In particular, when $n = 1$, we obtain
   \[
   (d^1_0 \lambda_1)(x \otimes y)
   \,=\, \lambda(x) \circ \lambda_1(y) - \lambda_1(xy) + \lambda_1(x) \circ \lambda(y)
   \,=\, 0.
   \]
Likewise, the condition \eqref{eq2:def} can be restated as
   \begin{equation}
   \label{eq:def2}
   \lambda_n(x)(ab) \,=\, \sum_{(x)}\sum_{i+j\,=\,n} \lambda_i(x_{(1)})(a) \cdot \lambda_j(x_{(2)})(b)
   \end{equation}
for all $n \geq 0$, $x \in H$, and $a, b \in A$.  By a simple rearrangement of terms, the case $n = 1$ then states that $d^1_1 \lambda_1 = 0$.  Therefore, we have $d^1 \lambda_1 = 0 \in \cF^2(A)$.  The last assertion about $\lambda_{k+1}$ is proved by essentially the same argument.

Now suppose that $\Lambdabar_t = \Phi^{-1}_t \Lambda_t \Phi_t$ for some deformation $\Lambdabar_t$ and formal automorphism $\Phi_t$.  Then the condition on the coefficients of $t$ is
   \[
   \lambdabar_1 \,=\, \lambda_1 + \lambda\circ \phi_1 - \phi_1 \circ \lambda,
   \]
i.e., $\lambdabar_1 - \lambda_1 = d^0 \phi_1$, a $1$-coboundary.  Therefore, the cohomology classes of $\lambda_1$ and $\lambdabar_1$ are the same.
\end{proof}

\subsection{Rigidity}
\label{subsec:rigidity}

The \emph{trivial deformation of $A$} is the deformation $\Lambda_t = \lambda$.  The $H$-module-algebra $A$ is said to be \emph{rigid} if and only if every deformation of $A$ is equivalent to the trivial deformation.

The following preliminary result is needed for the cohomological criterion for rigidity below.

\begin{prop}
\label{prop:rigidity}
Let $\Lambda_t = \lambda + \lambda_N t^N + O(t^{N+1})$ be a deformation of $A$ in which $\lambda_N = d^0 \phi$ for some $\phi \in \cF^0(A)$.  Suppose that $H^2_h(A,A) = 0$.  Then there exists a formal automorphism of $A$ of the form
   \begin{equation}
   \label{eq:aut}
   \Phi_t \,=\, \Id_A - \phi t^N + O(t^{N+1})
   \end{equation}
such that the deformation defined by
   \[
   \Lambdabar_t \,=\, \Phi^{-1}_t \Lambda_t \Phi_t
   \]
satisfies
   \[
   \lambdabar_i \,=\, 0
   \]
for $i = 1, \ldots, N$.
\end{prop}

\begin{proof}
The existence of a formal automorphism $\Phi_t$ of the form \eqref{eq:aut} is exactly \cite[Corollary 4.4]{yau}.  Since $\Phi_t \equiv \Id_A \pmod{t^N}$, we have
   \[
   \Lambdabar_t \,\equiv\, \Lambda_t \pmod{t^N}
   \]
and, therefore,
   \[
   \lambdabar_1 \,=\, \cdots \,=\, \lambdabar_{N-1} \,=\, 0.
   \]
To finish the proof, it suffices to consider the coefficient of $t^N$ in $\Phi^{-1}_t \Lambda_t \Phi_t$.  This coefficient is
   \[
   \lambdabar_N \,=\, \lambda_N + \phi \circ \lambda - \lambda \circ \phi
   \,=\, \lambda_N - d^0 \phi \,=\, 0,
   \]
as desired.
\end{proof}

Applying Theorem \ref{thm:inf} and Proposition \ref{prop:rigidity} repeatedly, we obtain the following cohomological criterion for the rigidity of a module-algebra.

\begin{cor}
\label{cor:rigidity}
If both $H^2_h(A,A)$ and $H^1(\cF^\bullet(A))$ are trivial, then $A$ is rigid.
\end{cor}


\section{Extending cocycles to deformations}
\label{sec:ext}

In view of Theorem \ref{thm:inf}, it is natural to ask the question:  Given a $1$-cocycle $\lambda_1 \in \cF^1(A)$, does there exist a deformation $\Lambda_t$ of $A$ whose infinitesimal is $\lambda_1$?  Following \cite{ger1}, if such a deformation exists, we say that $\lambda_1$ is \emph{integrable}.  The purpose of this section is to develop the obstruction theory for integrability of $1$-cocycles in $\cF^1(A)$.

\subsection{Deformations of finite order}
\label{subsec:finite}

Let $N$ be a positive integer.  A polynomial $\Lambda_t = \sum_{i=0}^N \lambda_i t^i$ with $\lambda_0 = \lambda$ and each $\lambda_i \in \cF^1(A)$ is said to be a \emph{deformation of order $N$} if and only if it satisfies the definition of a deformation modulo $t^{N+1}$, i.e., \eqref{eq:def1} and \eqref{eq:def2} for $n \leq N$.  Such a deformation of order $N$ is said to \emph{extend to order $N + 1$} if and only if there exists a $1$-cochain $\lambda_{N+1} \in \cF^1(A)$ such that the polynomial
   \begin{equation}
   \label{eq:extended poly}
   \Lambdatilde_t = \Lambda_t + \lambda_{N+1}t^{N+1} = \sum_{i=0}^{N+1} \lambda_i t^i
   \end{equation}
is a deformation of order $N + 1$.  In this case, we say that $\Lambdatilde_t$ is an \emph{order $N + 1$ extension of $\Lambda_t$}.

It is easy to see that a $1$-cochain $\lambda_1 \in \cF^1(A)$ is a $1$-cocycle if and only if the linear polynomial $\lambda + \lambda_1 t$ is a deformation of order $1$.  Therefore, in order to find the obstructions to integrating $\lambda_1$, it suffices to find the obstruction to extending a deformation of order $N \geq 1$ to one of order $N + 1$.

\subsection{Obstruction}
\label{subsec:ob}

Let, then, $\Lambda_t = \sum_{i=0}^N \lambda_i t^i$ be a deformation of order $N \geq 1$.  As explained in the proof of Theorem \ref{thm:inf}, the conditions \eqref{eq:def1} and \eqref{eq:def2} for $n = 1$ are equivalent to $\lambda_1 \in \cF^1(A)$ being a $1$-cocycle.  For each $m = 2, \ldots, N$, the conditions \eqref{eq:def1} and \eqref{eq:def2} for $n = m$, by a simple rearrangement of terms, are equivalent to
   \begin{equation}
   \label{eq:ob m}
   d^1 \lambda_m \,=\, - \sum_{i=1}^{m-1} \lambda_i \cup \lambda_{m-i},
   \end{equation}
where the $\cup$-product was introduced in \eqref{eq:cup}.

Let $\lambda_{N+1} \in \cF^1(A)$ be an arbitrary $1$-cochain and set $\Lambdatilde_t = \Lambda_t + \lambda_{N+1}t^{N+1}$.  Then $\Lambdatilde_t$ is a deformation of order $N + 1$ if and only if it satisfies \eqref{eq:def1} and \eqref{eq:def2} for $n = N + 1$.  As in the previous paragraph, this is equivalent to
   \begin{equation}
   \label{eq:ob N+1}
   d^1 \lambda_{N+1} \,=\, - \sum_{i=1}^N \lambda_i \cup \lambda_{N+1-i}.
   \end{equation}
Consider the $2$-cochain
   \begin{equation}
   \label{eq:ob}
   \Ob \,=\, \sum_{i=1}^N \lambda_i \cup \lambda_{N+1-i}
   \end{equation}
in $\cF^2(A)$ defined by $\lambda_1, \ldots, \lambda_N$.

\begin{lemma}
\label{lem:ob}
The class $\Ob \in \cF^2(A)$ is a $2$-cocycle.
\end{lemma}

\begin{proof}
This is similar to \cite[Proposition 2]{ger1}.  In fact, using $d^1 \lambda_1 = 0$, \eqref{eq:ob m}, and Proposition \ref{prop:cup}, we have:
   \[
   \begin{split}
   d^2 \Ob
   & \,=\, \sum_{i=1}^N d^2(\lambda_i \cup \lambda_{N+1-i}) \\
   & \,=\, \sum_{i=1}^N \left\lbrace (d^1\lambda_i) \cup \lambda_{N+1-i} - \lambda_i \cup (d^1\lambda_{N+1-i}) \right\rbrace \\
   & \,=\, \sum_{i=2}^N \left\lbrace - \sum_{j=1}^{i-1} \lambda_j \cup \lambda_{i-j} \right\rbrace \cup \lambda_{N+1-i} - \sum_{i=1}^N \lambda_i \cup \left\lbrace - \sum_{j=1}^{N-i} \lambda_j \cup \lambda_{N+1-i-j} \right\rbrace \\
   & \,=\, - \sumprime \lambda_a \cup \lambda_b \cup \lambda_c + \sumprime \lambda_a \cup \lambda_b \cup \lambda_c \\
   & \,=\, 0.
   \end{split}
   \]
Here $\sumprime$ is the sum over all integers $a, b, c > 0$ with $a + b + c = N + 1$.
\end{proof}

Combining \eqref{eq:ob N+1}, \eqref{eq:ob}, and Lemma \ref{lem:ob}, we obtain the desired obstruction for extending an order $N$ deformation.

\begin{thm}
\label{thm:ext}
The deformation $\Lambda_t$ of order $N$ extends to order $N + 1$ if and only if the $2$-cocycle $-\Ob$ is a $2$-coboundary.
\end{thm}

Since the obstruction is always a class in $H^2(\cF^\bullet(A))$, we obtain the following cohomological criterion for integrability.

\begin{cor}
\label{cor:ext}
If $H^2(\cF^\bullet(A))$ is trivial, then every $1$-cocycle in $\cF^1(A)$ is integrable.
\end{cor}


\section{Deformation cohomology for module-coalgebras and comodule-(co)algebras}
\label{sec:def'}

The purpose of this final section is to describe the deformation DGAs and the corresponding deformation results for module-coalgebras, comodule-algebras, and comodule-coalgebras.  As before, the deformation is taken with respect to the module (or comodule) action, leaving the algebra (or coalgebra) structure unaltered.  In each case, once the correct deformation DGA is set up, the statements of results and their arguments are formally similar to the module-algebra case above.  Therefore, we will describe the constructions and statements of results and omit the arguments.  To avoid too much repetitions, we will concentrate on the case of comodule-coalgebras.  At the end of the section, we will indicate what modifications are needed for the cases of module-coalgebras and comodule-algebras.

\subsection{Comodule-coalgebras}
\label{subsec:cc}

We still denote by $(H, \mu_H, \Delta_H)$ a bialgebra.  Let $(A, \Delta_A)$ be a coalgebra.  An \emph{$H$-comodule-coalgebra} structure on $A$ consists of an $H$-comodule structure $\rho \colon A \to H \otimes A$ on $A$ such that the map $\Delta_A \colon A \to A \otimes A$ is an $H$-comodule morphism, i.e.,
   \begin{equation}
   \label{eq:cc}
   (\Id_H \otimes \Delta_A) \circ \rho
   \,=\, (\mu_H \otimes \Id_{A^{\otimes 2}}) \circ (\Id_H \otimes \tau \otimes \Id_A) \circ \rho^{\otimes 2} \circ \Delta_A.
   \end{equation}
Here $\tau \colon A \otimes H \cong H \otimes A$ is the twist isomorphism.

Until otherwise indicated, $A$ will denote an $H$-comodule-coalgebra with structure map $\rho$.

\subsection{Deformations of comodule-coalgebras}
\label{subsec:def cc}

A \emph{deformation} of $A$ is a power series $R_t = \sum_{i=0}^\infty \rho_i t^i$ with $\rho_0 = \rho$ and each $\rho_i \in \Hom(A, H \otimes A)$, satisfying the following two conditions:
   \begin{subequations}
   \label{eq:def cc}
   \begin{align}
   (\Id_H \otimes R_t) \circ R_t & \,=\, (\Delta_H \otimes \Id_A) \circ R_t, \label{eq1:def cc} \\
   (\Id_H \otimes \Delta_A) \circ R_t & \,=\, (\mu_H \otimes \Id_{A^{\otimes 2}}) \circ (\Id_H \otimes \tau \otimes \Id_A) \circ R_t^{\otimes 2} \circ \Delta_A. \label{eq2:def cc}
   \end{align}
   \end{subequations}
As before, maps are extended linearly to include modules of power series wherever appropriate.  The linear coefficient $\rho_1$ is called the \emph{infinitesimal} of $R_t$.

A \emph{formal automorphism} of $A$ is a power series $\Phi_t = \sum_{i=0}^\infty \phi_i t^i$ with $\phi_0 = \Id_A$ and each $\phi_i \in \End(A)$ that is comultiplicative, i.e.,
   \[
   \Delta_A \circ \Phi_t \,=\, \Phi_t^{\otimes 2} \circ \Delta_A.
   \]
Two deformations $R_t$ and $\Rbar_t$ of $A$ are \emph{equivalent} if and only if there exists a formal automorphism of $A$ such that
   \begin{equation}
   \label{eq:equiv cc}
   \Rbar_t \,=\, (\Id_H \otimes \Phi_t^{-1}) \circ R_t \circ \Phi_t.
   \end{equation}

\subsection{Deformation complex for a comodule-coalgebra}
\label{subsec:DGA cc}

It is the cochain complex $(\cF^\bullet_{cc}(A), d_{cc})$ defined as follows:
   \[
   \begin{split}
   \cF^0_{cc}(A) & \,=\, \Coder(A), \\
   \cF^1_{cc}(A) & \,=\, \Hom(A, H \otimes A), \\
   \cF^n_{cc}(A) & \,=\, \cF^n_{cc,0}(A) \oplus \cF^n_{cc,1}(A) \quad (n \geq 2),
   \end{split}
   \]
where
   \[
   \begin{split}
   \cF^n_{cc,0}(A) & \,=\, \Hom(A, H^{\otimes n} \otimes A), \\
   \cF^n_{cc,1}(A) & \,=\, \Hom(A, H \otimes A^{\otimes n}).
   \end{split}
   \]
Now we define the differentials.
   \[
   \begin{split}
   d^0_{cc}\phi & \,=\, \rho \circ \phi - (\Id_H \otimes \phi) \circ \rho, \\
   d^n_{cc} & \,=\, (d^n_{cc,0}; d^n_{cc,1}) \quad (n \geq 1).
   \end{split}
   \]
The component maps $d^n_{cc,i} \colon \cF^n_{cc,i}(A) \to \cF^{n+1}_{cc,i}(A)$ $(i = 0, 1)$ are defined by:
   \[
   \begin{split}
   (d^n_{cc,0}\varphi_0) & \,=\, (\Id_{H^{\otimes n}} \otimes \rho) \circ \varphi_0
   + \sum_{i=1}^n (-1)^i (\Id_{H^{\otimes (n-i)}} \otimes \Delta_H \otimes \Id_{H^{\otimes (i-1)} \otimes A}) \circ \varphi_0 \\
   & \relphantom{} \relphantom{} + (-1)^{n+1}(\Id_H \otimes \varphi_0) \circ \rho  \\
   (d^n_{cc,1}\varphi_1) & \,=\, (\mu_H \otimes \Id_{A^{\otimes (n+1)}}) \circ (\Id_H \otimes \tau \otimes \Id_{A^{\otimes n}}) \circ (\rho \otimes \varphi_1) \circ \Delta_A \\
   & \relphantom{} \relphantom{} + \sum_{i=1}^n (-1)^i (\Id_{H \otimes A^{\otimes (i-1)}} \otimes \Delta_A \otimes \Id_{A^{\otimes (n-i)}}) \circ \varphi_1 \\
   & \relphantom{} \relphantom{} + (-1)^{n+1} (\mu_H \otimes \Id_{A^{\otimes (n+1)}}) \circ (\Id_H \otimes \tau_n \otimes \Id_A) \circ (\varphi_1 \otimes \rho) \circ \Delta_A.
   \end{split}
   \]
In the last line, the map $\tau_n \colon A^{\otimes n} \otimes H \cong H \otimes A^{\otimes n}$ is the twist isomorphism
   \[
   \tau_n(\mathbf{a} \otimes x) \,=\, x \otimes \mathbf{a}.
   \]

\subsection{Cup-product on $\cF^\bullet_{cc}(A)$}
\label{subsec:cup cc}

There is a $\cup$-product
   \begin{equation}
   \label{eq:cup cc}
   - \cup - \colon \cF^m_{cc}(A) \otimes \cF^n_{cc}(A) \to \cF^{m+n}_{cc}(A) \quad (m, n > 0)
   \end{equation}
that is defined as follows.  For $f \in \cF^m_{cc}(A)$ and $g \in \cF^n_{cc}(A)$, the components of $f \cup g$ are:
   \[
   \begin{split}
   (f \cup g)_0 & \,=\, (\Id_{H^{\otimes n}} \otimes f_0) \circ g_0, \\
   (f \cup g)_1 & \,=\, (\mu_H \otimes \Id_{A^{\otimes (m+n)}}) \circ (\Id_H \otimes \tau_m \otimes \Id_{A^{\otimes n}}) \circ (f_1 \otimes g_1) \circ \Delta_A.
   \end{split}
   \]

\begin{thm}
\label{thm:cochain cc}
$(\cF^\bullet_{cc}(A), d_{cc})$ is a cochain complex.  Moreover, $(\cF^{> 0}_{cc}(A), d_{cc}, \cup)$ is a DGA, and $(H^{> 0}(\cF^\bullet_{cc}(A)), \cup)$ is a graded algebra.
\end{thm}

\subsection{Hochschild coalgebra cohomology}

The rigidity result for a comodule-coalgebra (and also a module-coalgebra) uses Hochschild coalgebra cohomology \cite{jonah,pw}, which we now recall.  The \emph{Hochschild coalgebra cohomology} of a coalgebra $A$ with coefficients in an $A$-bicomodule $M$ (with left $A$-coaction $\psi_l$ and right $A$-coaction $\psi_r$) is defined as follows.  For $n \geq 1$, the module of $n$-cochains is defined to be
   \[
   C_c^n(M,A) = \Hom(M, A^{\otimes n}),
   \]
with differential
   \[
   \begin{split}
   \delta_c \sigma
   & = (\Id_A \otimes \sigma) \circ \psi_l  \,+\, \sum_{i=1}^n (-1)^i \left(\Id_{A^{\otimes(i-1)}} \otimes \Delta \otimes \Id_{A^{\otimes(n-i)}}\right) \circ \sigma \\
   & \relphantom{} + (-1)^{n+1} (\sigma \otimes \Id_A) \circ \psi_r
   \end{split}
   \]
for $\sigma \in C^n_c(M,A)$.  Set $C^0_c(M,A) \equiv 0$.  The cohomology of the cochain complex $(C^\bullet_c(M,A), \delta_c)$ is denoted by $H^\bullet_c(M,A)$.

For example, we can consider $A$ as an $A$-bicomodule with coaction maps $\psi_l = \psi_r = \Delta_A$.

\subsection{Results about deformations of comodule-coalgebras}
\label{subsec:def result}

The \emph{trivial deformation} of $A$ is the deformation $R_t = \rho$.  Rigidity is defined as in the module-algebra case.

\begin{thm}
\label{thm:def cc}
Let $R_t = \sum_{i=0}^\infty \rho_i t^i$ be a deformation of $A$ as an $H$-comodule-coalgebra.  Then the infinitesimal $\rho_1$ is a $1$-cocycle in $\cF^1_{cc}(A)$ whose cohomology class is determined by the equivalence class of $R_t$.  Moreover:
   \begin{enumerate}
   \item If $H^2_c(A,A)$ and $H^1(\cF^\bullet_{cc}(A))$ are both trivial, then $A$ is rigid.
   \item If $H^2(\cF^\bullet_{cc}(A))$ is trivial, then every $1$-cocycle in $\cF^1_{cc}(A)$ is the infinitesimal of some deformation of $A$.
   \end{enumerate}
\end{thm}

\subsection{Module-coalgebras}
\label{subsec:mc}

Next we consider the case of module-coalgebras.  An \emph{$H$-module-coalgebra} structure on the coalgebra $A$ consists of an $H$-module structure $\lambda \in \Hom(H, \End(A))$ on $A$ such that the map $\Delta_A \colon A \to A \otimes A$ becomes an $H$-module morphism, i.e.,
   \[
   \Delta_A(\lambda(x)(a)) \,=\, \sum_{(a)(x)} \lambda(x_{(1)})(a_{(1)}) \otimes \lambda(x_{(2)})(a_{(2)})
   \]
for all $x \in H$ and $a \in A$.

A \emph{deformation} $\Lambda_t = \sum_{i=0}^\infty \lambda_i t^i$ of $A$ is defined as in the module-algebra case, except that \eqref{eq2:def} is replaced by
   \begin{equation}
   \label{eq2:def mc}
   \Delta_A(\Lambda_t(x)(a)) \,=\, \sum_{(a)(x)} \Lambda_t(x_{(1)})(a_{(1)}) \otimes \Lambda_t(x_{(2)})(a_{(2)}).
   \end{equation}

A \emph{formal automorphism} of $A$ is defined as in the comodule-coalgebra case (\S \ref{subsec:def cc}), and \emph{equivalence} is defined as in the module-algebra case (\S \ref{subsec:equivalence}).

The deformation DGA $(\cF^\bullet_{mc}(A), d_{mc}, \cup)$ of $A$ (as an $H$-module-coalgebra with structure map $\lambda$) has the same general form as in the cases of module-algebra and comodule-coalgebra:
   \[
   \allowdisplaybreaks
   \begin{split}
   \cF^0_{mc}(A) & \,=\, \Coder(A), \\
   \cF^1_{mc}(A) & \,=\, \Hom(H, \End(A)), \\
   \cF^n_{mc,0}(A) & \,=\, \cF^n_0(A) \quad (n \geq 2), \\
   \cF^n_{mc,1}(A) & \,=\, \Hom(H, \Hom(A, A^{\otimes n})) \quad (n \geq 2).
   \end{split}
   \]
The differentials are defined as follows:
   \[
   \begin{split}
   d^0_{mc}\varphi & \,=\, \lambda \circ \phi - \phi \circ \lambda, \\
   d^n_{mc,0} & \,=\, d^n_0 \quad (n \geq 1), \\
   (d^n_{mc,1}\varphi_1)(x)(a) & \,=\, \sum_{(a)(x)} \lambda(x_{(1)})(a_{(1)}) \otimes \varphi_1(x_{(2)})(a_{(2)}) \\
   & \relphantom{} \relphantom{} + \sum_{i=1}^n (-1)^i (\Id_{A^{\otimes (i-1)}} \otimes \Delta_A \otimes \Id_{A^{\otimes (n-i)}})(\varphi_1(x)(a)) \\
   & \relphantom{} \relphantom{} + (-1)^{n+1} \sum_{(a)(x)} \varphi_1(x_{(1)})(a_{(1)}) \otimes \lambda(x_{(2)})(a_{(2)}) \quad (n \geq 1).
   \end{split}
   \]
The $\cup$-product is defined by:
   \[
   \begin{split}
   (f \cup g)_0 & \,=\, f_0 \cup g_0 \quad (\text{in } C^{> 0}(H, \End(A))), \\
   (f \cup g)_1(x)(a) & \,=\, \sum_{(a)(x)} f_1(x_{(1)})(a_{(1)}) \otimes g_1(x_{(2)})(a_{(2)}).
   \end{split}
   \]
With these definitions, we have the module-coalgebra analogue of Theorem \ref{thm:cochain cc} and Theorem \ref{thm:def cc}, where $(\cF^\bullet_{cc}(A), d_{cc})$ is replaced by $(\cF^\bullet_{mc}(A), d_{mc})$.

\subsection{Comodule-algebras}
\label{subsec:ca}

Finally, we consider the case of comodule-algebras.  Here let $(A, \mu_A)$ be an algebra.  An \emph{$H$-comodule-algebra} structure on $A$ consists of an $H$-comodule structure $\rho \colon A \to H \otimes A$ on $A$ such that the map $\mu_A \colon A \otimes A \to A$ becomes an $H$-comodule morphism, i.e.,
   \[
   \rho \circ \mu_A \,=\, (\mu_H \otimes \mu_A) \circ (\Id_H \otimes \tau \otimes \Id_A) \circ \rho^{\otimes 2}.
   \]

A \emph{deformation} of $A$ (as an $H$-comodule-algebra with structure map $\rho$) is defined as in the case of comodule-coalgebras (\S \ref{subsec:def cc}), except that the condition \eqref{eq2:def cc} is replaced by
   \[
   R_t \circ \mu_A \,=\, (\mu_H \otimes \mu_A) \circ (\Id_H \otimes \tau \otimes \Id_A) \circ R_t^{\otimes 2}.
   \]

A \emph{formal automorphism} of $A$ is defined as in the module-algebra case (\S \ref{subsec:equivalence}), and \emph{equivalence} is defined as in the case of comodule-coalgebras \eqref{eq:equiv cc}.

The deformation DGA $(\cF^\bullet_{ca}(A), d_{ca}, \cup)$ of $A$ as an $H$-comodule-algebra takes the following form:
   \[
   \allowdisplaybreaks
   \begin{split}
   (\cF^0_{ca}(A), d^0_{ca}) & \,=\, (\Der(A), d^0_{cc}) \\
   \cF^1_{ca}(A) & \,=\, \Hom(A, H \otimes A) \\
   \cF^n_{ca}(A) & \,=\, \cF^n_{ca,0}(A) \oplus \cF^n_{ca,1}(A) \quad (n \geq 2) \\
   (\cF^n_{ca,0}(A), d^n_{ca,0}) & \,=\, (\cF^n_{cc,0}(A), d^n_{cc,0}) \quad (n \geq 1) \\
   \cF^n_{ca,1}(A) & \,=\, \Hom(A^{\otimes n}, H \otimes A) \quad (n \geq 2) \\
   d^n_{ca,1}\varphi_1 & \,=\, (\mu_H \otimes \mu_A) \circ (\Id_H \otimes \tau \otimes \Id_A) \circ (\rho \otimes \varphi_1) \\
   & \relphantom{} \relphantom{} + \sum_{i=1}^n (-1)^i \varphi_1 \circ (\Id_{A^{\otimes (i-1)}} \otimes \mu_A \otimes \Id_{A^{\otimes (n-i)}}) \\
   & \relphantom{} \relphantom{} + (-1)^{n+1} (\mu_H \otimes \mu_A) \circ (\Id_H \otimes \tau \otimes \Id_A) \circ (\varphi_1 \otimes \rho) \quad (n \geq 1).
   \end{split}
   \]
The $\cup$-product is given by:
   \[
   \begin{split}
   (f \cup g)_0 & \,=\, (\Id_{H^{\otimes n}} \otimes f_0) \circ g_0 \\
   (f \cup g)_1 & \,=\, (\mu_H \otimes \mu_A) \circ (\Id_H \otimes \tau \otimes \Id_A) \circ (f_1 \otimes g_1).
   \end{split}
   \]
With these definitions, we have the comodule-algebra analogue of Theorem \ref{thm:cochain cc} and Theorem \ref{thm:def cc}, where $(\cF^\bullet_{cc}(A), d_{cc})$ and $H^2_c(A,A)$ are replaced by $(\cF^\bullet_{ca}(A), d_{ca})$ and $H^2_h(A,A)$, respectively.


\end{document}